\documentclass{amsart}
\usepackage{lastpage}
\usepackage{yhmath,verbatim}
\usepackage[all]{xy}

\newcommand{\bdis}{\begin{displaymath}}
\newcommand{\edis}{\end{displaymath}}
\newcommand{\be}{\begin{equation}}
\newcommand{\ee}{\end{equation}}
\newcommand{\mbb}{\mathbb}
\newcommand{\mcal}{\mathcal}

\newcommand{\vp}{\varphi}

\newcommand{\zf}{\zeta\left(\frac{1}{2}+it\right)}

\DeclareMathOperator{\im}{Im}

\theoremstyle{definition}

\theoremstyle{remark}
\newtheorem{remark}[]{Remark}

\newtheorem*{mydef11}{{\bf Theorem 1}}

\newtheorem*{mydef21}{{\bf Definition 1}}

\newtheorem*{mydef41}{{\bf Corollary 1}}

\newtheorem*{mydef51}{{\bf Lemma 1}}

\newtheorem*{mydef52}{{\bf Lemma 2}}

\newtheorem*{mydef53}{{\bf Lemma 3}}

\newtheorem*{mydef54}{{\bf Lemma 4}}

\newtheorem*{mydef55}{{\bf Lemma 5}}

\newtheorem*{mydef56}{{\bf Lemma 6}}

\newtheorem*{mydef57}{{\bf Lemma 7}}

\newtheorem*{mydef81}{{\bf Property 1}}

\newtheorem*{mydef82}{{\bf Property 2}}

\numberwithin{equation}{section}

\begin{document}

\title{Jacob's ladder as generator of new class of iterated $L_2$-orthogonal systems and their dependence on the Riemann's function $\zf$} 

\author{Jan Moser}

\address{Department of Mathematical Analysis and Numerical Mathematics, Comenius University, Mlynska Dolina M105, 842 48 Bratislava, SLOVAKIA}

\email{jan.mozer@fmph.uniba.sk}


\begin{abstract}
In this paper new classes of $L_2$-orthogonal functions are constructed as iterated $L_2$-orthogonal systems. In order to do this we use the theory of the Riemann's zeta-function as well as our theory of Jacob's ladders. The main result is new one in the theory of the Riemann's zeta-function and simultaneously in the theory of $L_2$-orthogonal systems.  \\
\vspace{0.2cm} \ 

\noindent 
DEDICATED TO THE MEMORY OF FOURIER'S EGYPTIAN ANABASE 

\end{abstract}
\maketitle

\section{Main result} 

\subsection{}  

Let us remind the following notions: 
\begin{itemize}
	\item[(a)] Jacob's ladder $\vp_1(t)$, 
	\item[(b)] the function 
	\be \label{1.1} 
	\begin{split}
	& \tilde{Z}^2(t)=\frac{{\rm d}\vp_1(t)}{{\rm d}t}=\frac{1}{\omega(t)}\left|\zf\right|^2, \\ 
	& \omega(t)=\left\{1+\mcal{O}\left(\frac{\ln\ln t}{\ln t}\right)\right\}\ln t,\ t\to\infty, 
	\end{split}
	\ee 
	\item[(c)] the (direct) iterations of the Jacob's ladder 
	\bdis 
	\vp_1^0(t)=t,\ \vp_1^1(t)=\vp_1(t),\ \vp_1^2(t)=\vp_1(\vp_1(t)),\ \dots ,\vp_1^k(t)=\vp_1(\vp_1^{k-1}(t)) 
	\edis  
	for every fixed $k\in\mbb{N}$, 
	\item[(d)] the reverse iterations (by means of $\vp_1(t)$) 
	\bdis 
	\begin{split}
	& [\overset{0}{T},\overset{0}{\wideparen{T+U}}],\ [\overset{1}{T},\overset{1}{\wideparen{T+U}}],\ \dots,\ [\overset{k}{T},\overset{k}{\wideparen{T+U}}], \\ 
	& U=o\left(\frac{T}{\ln T}\right),\ T\to\infty, 
	\end{split}
	\edis  
	of the basic segment 
	\bdis 
	[T,T+U]=[\overset{0}{T},\overset{0}{\wideparen{T+U}}], 
	\edis  
	where 
	\be \label{1.2} 
	[\overset{0}{T},\overset{0}{\wideparen{T+U}}]\prec [\overset{1}{T},\overset{1}{\wideparen{T+U}}]\prec \dots \prec 
	[\overset{k}{T},\overset{k}{\wideparen{T+U}}],  
	\ee  
	that we have introduced into the theory of the Riemann's zeta-function in our papers \cite{1} -- \cite{4}. 
\end{itemize} 

\subsection{} 

The following theorem is the main result of this paper. 

\begin{mydef11}
For every fixed $L_2$-orthogonal system 
\be \label{1.3} 
\{f_n(t)\}_{n=0}^\infty,\ t\in [a,a+2l],\ a\in\mbb{R},\ l\in\mbb{R}^{+} 
\ee  
and for every fixed $k\in\mbb{N}$ there is the set of $k$ new iterated $L_2$-orthogonal systems 
\be \label{1.4} 
\{f_n^p(t)\}_{n=0}^\infty,\ t\in [a,a+2l],\ p=1,2,\dots,k , 
\ee  
where 
\be \label{1.5} 
\begin{split}
& f_n^p(t)=f_n\left(\vp_1^p\left(\frac{\overset{p}{\wideparen{T+2l}}-\overset{p}{T}}{2l}(t-a)+\overset{p}{T}\right)-T+a\right)\times \\ 
& \times \prod_{r=0}^{p-1}\left|\tilde{Z}\left(\vp_1^r\left(\frac{\overset{p}{\wideparen{T+2l}}-\overset{p}{T}}{2l}(t-a)+\overset{p}{T}\right)\right)\right|,\ 
\end{split}
\ee 
for all sufficiently big $T>0$. 
\end{mydef11}  

\subsection{} 

For example, by Theorem 1, we can assign the following set of $k$ new iterated orthogonal systems 
\bdis 
\begin{split}
	& P_n^p(t)=P_n\left(\vp_1^p\left(\frac{\overset{p}{\wideparen{T+2}}-\overset{p}{T}}{2l}(t+1)+\overset{p}{T}\right)-T-1\right)\times \\ 
	& \times \prod_{r=0}^{p-1}\left|\tilde{Z}\left(\vp_1^r\left(\frac{\overset{p}{\wideparen{T+2}}-\overset{p}{T}}{2l}(t+1)+\overset{p}{T}\right)\right)\right|, 
\end{split}
\edis 
for all sufficiently big $T$ 
to the classical Legendre's orthogonal system 
\bdis 
\{P_n(t)\}_{n=0}^\infty,\ t\in[-1,1];\ a=-1,\ l=1. 
\edis  
For example 
\bdis 
\begin{split}
	& P_n^1(t)=P_n\left(\vp_1\left(\frac{\overset{1}{\wideparen{T+2}}-\overset{1}{T}}{2l}(t+1)+\overset{1}{T}\right)-T-1\right)\times \\ 
	& \times  \left|\tilde{Z}\left(\vp_1\left(\frac{\overset{1}{\wideparen{T+2}}-\overset{1}{T}}{2l}(t+1)+\overset{1}{T}\right)\right)\right|, \ t\in[-1,1]. 
\end{split}
\edis 

\subsection{} 

Restating Theorem 1 we have the following 

\begin{mydef41} 
For $L_2$-orthogonal system (\ref{1.3}) there is the set of $k$ new iterated $L_2$-orthogonal systems 
\bdis 
\begin{split}
	& \left\{f_n\left(\vp_1^p\left(\frac{\overset{p}{\wideparen{T+2l}}-\overset{p}{T}}{2l}(t-a)+\overset{p}{T}\right)-T+a\right)\right\}_{n=0}^\infty , \\ 
	& t\in [a,a+2l],\ p=1,\dots,k, 
\end{split}
\edis 
with weights 
\bdis 
\prod_{r=0}^{p-1}\tilde{Z}^2\left(\vp_1^r\left(\frac{\overset{p}{\wideparen{T+2}}-\overset{p}{T}}{2l}(t+1)+\overset{p}{T}\right)\right), 
\edis  
and this last is (see (\ref{1.1})) 
\bdis 
\sim \frac{1}{\ln^pT}\prod_{r=0}^{p-1}\left|
\zeta\left(\frac 12+i\vp_1^r\left(\frac{\overset{p}{\wideparen{T+2l}}-\overset{p}{T}}{2l}(t-a)+\overset{p}{T}\right)\right)\right|^2,\ T\to\infty. 
\edis 
\end{mydef41}  

\subsection{} 

Now we give some remarks. 

\begin{remark}
Theorem 1 represents completely new result in the theory of the Riemann's zeta-function and simultaneously in the the theory of $L_2$-orthogonal systems. 
\end{remark} 

\begin{remark}
The last row \emph{for all sufficiently big $T>0$} in the Theorem 1 gives the continuum set of possibilities for construction of new $k$-tuples of iterated $L_2$-orthogonal systems for every fixed system (\ref{1.3}). 
\end{remark} 

\begin{remark}
Dependence of iterated $L_2$-orthogonal systems (\ref{1.4}) on the Riemann's zeta-function $\zf$ is evident one, see (\ref{1.1}), (\ref{1.5}). 
\end{remark} 

\begin{remark}
This paper is the continuation of 54 papers concerning Jacob's ladders. These can be found in arXiv [math.CA] starting with the paper \cite{1}. 
\end{remark}

\section{Jacob's ladders} 

\subsection{} 

Let us remind the following non-linear integral equation 
\be \label{2.1} 
\int_0^{\mu[x(T)]}Z^2(t)e^{-\frac{2}{x(T)}t}{\rm d}t=\int_0^T Z^2(t){\rm d}t 
\ee  
we have introduced in the paper \cite{1}, where 
\be \label{2.2} 
\begin{split}
& Z(t)=e^{i\vartheta(t)}\zf , \\ 
& \vartheta(t)=-\frac t2\ln\pi+\im\left\{\ln\Gamma\left(\frac 14+i\frac t2\right)\right\}, 
\end{split}
\ee 
and the class of functions $\{\mu\}$ specified as 
\bdis 
\mu\in C^\infty([y_0,+\infty)) 
\edis 
monotonically increasing and unbounded from above and obeying the inequality 
\be \label{2.3} 
\mu(y)\geq 7y\ln y. 
\ee 

\subsection{} 

The following statement holds (see \cite {1}). 

\begin{mydef51}
For any $\mu\in\{\mu\}$ there is exactly one solution to the integral equation (\ref{2.1}): 
\bdis 
\begin{split}
& \vp(T)=\vp(T,\mu),\ T\in [T_0,+\infty),\ T_0=T_0[\vp]>0, \\ 
& \vp(T)\xrightarrow{T\to\infty}\infty. 
\end{split}
\edis 
\end{mydef51} 

Let the symbol $\{\vp\}$ denote the system of these solutions. The function $\vp(T)$ is related to the zeroes of the Riemann's zeta-function on the critical line by the following way. Let $t=\gamma$ be such a zero of 
\bdis 
\zf 
\edis 
and of the order $n(\gamma)$, where 
\bdis 
n(\gamma)=\mcal{O}(\ln \gamma),\ \gamma\to\infty, 
\edis   
of course. Then the following holds true. 

\begin{remark}
The points 
\bdis 
[\gamma,\vp(\gamma)],\ \gamma>T_0
\edis 
(and only these points) are inflection points with the horizontal tangent. In more detail, the following system of equations holds true: 
\bdis 
\vp'(\gamma)=\vp''(\gamma)=\dots=\vp^{(2n)}(\gamma)=0,\ \vp^{(2n+1)}(\gamma)\not=0,\ n=n(\gamma). 
\edis 
\end{remark} 

\begin{mydef21}
With respect the above mentioned property, an element $\vp\in\{\vp\}$ is called Jacob's ladder leading to $[+\infty,+\infty]$. The rungs of that ladder are segments of the curve $\vp$ lying in the neighbourhoods of the points 
\bdis 
[\gamma,\vp(\gamma)],\ \gamma>T_0[\vp].  
\edis 
\end{mydef21} 

\begin{remark}
We call $\vp(T)$ as Jacob's ladder in analogy with Jacob's dream in Chumash, Bereishis, 28:12. 
\end{remark} 

\begin{remark}
Finally, the composite function $g[\vp(T)]$ is also called Jacob's ladder for any function $g$ that is increasing, $C^\infty$ on $[y_0,+\infty)$ and unbounded from above. For example, the function 
\be \label{2.4} 
\vp_1(T)=\frac 12\vp(T) 
\ee  
as composition of $g(y)=\frac 12 y$, $y\geq y_0$, $y=\vp(T)$,\ $T\geq T_0[\vp]=T_0[\vp_1]$, $g'_y=\frac 12>0$ is the Jacob's ladder. 
\end{remark} 

\section{Basic property of Jacob's ladders: existence of almost exact expressions for the classical Hardy-Littlewood integral (1918)}

\subsection{} 

Let us remind the Hardy-Littlewood integral 
\be \label{3.1} 
\int_0^T\left|\zf\right|^2{\rm d}t 
\ee  
can be expressed as follows: 
\be \label{3.2} 
\int_0^T\left|\zf\right|^2{\rm d}t = T\ln T+(2c-1-\ln 2\pi)T+R(T), 
\ee  
with, for example, Ingham's error term 
\be \label{3.3} 
R(T)=\mcal{O}(T^{\frac 12}\ln T)=\mcal{O}(T^{\frac 12+\delta}),\ \delta>0,\ T\to\infty 
\ee  
for arbitrarily small $\delta$. 

Next, by Good's $\Omega$-theorem (1977) we have that 
\be \label{3.4} 
R(T)=\Omega(T^{\frac 14}),\ T\to\infty. 
\ee  

\begin{remark}
Let 
\be \label{3.5} 
R_a(T)=\mcal{O}(T^{\frac 14+a}),\ a\in [\delta,\frac 14+\delta],\ T\to\infty. 
\ee 
Then, by (\ref{3.4}), it is true that for every valid estimate of type (\ref{3.5}) one obtains: 
\be \label{3.6} 
\limsup_{T\to\infty}|R_a(T)|=+\infty. 
\ee  
In other words, every expression of the type (\ref{3.2}) possesses an unbounded error term at infinity. 
\end{remark} 

\subsection{} 

Under the circumstances (\ref{3.2}) and (\ref{3.6}) we have proved in \cite{1} that the Hardy-Littlewood integral (\ref{3.1}) has an infinite set of other completely new and almost exact representations expressed by the following: 

\begin{mydef81}
\be \label{3.7} 
\begin{split} 
& \int_0^T\left|\zf\right|^2{\rm d}t= \vp_1(T)\ln\{\vp_1(T)\}+ \\ 
& + (c-\ln 2\pi)\vp_1(T)+c_0+\mcal{O}\left(\frac{\ln T}{T}\right),\ T\to\infty, 
\end{split}
\ee  
(comp. (\ref{2.4})) with the following 
\be \label{3.8} 
\lim_{T\to\infty}\tilde{R}(T)=\lim_{T\to\infty}\mcal{O}\left(\frac{\ln T}{T}\right)=0, 
\ee  
where $c$ is the Euler's constant and $c_0$ is the constant from the Titchmarsh-Kober-Atkinson formula. 
\end{mydef81} 

\begin{remark}
Comparison of (\ref{3.6}) and (\ref{3.8}) completely characterizes the level of exactness of our representation (\ref{3.7}) of the Hardy-Littlewood integral (\ref{3.1}). 
\end{remark} 

\section{Asymptotic relation between Jacob's ladder and the prime-counting function} 

\subsection{} 

Further, in the paper \cite{1}, (6.2), we have obtained the following formula. 

\begin{mydef52}
\be\label{4.1} 
T-\vp_1(T)\sim (1-c)\pi(T);\ \pi(T)\sim \frac{T}{\ln T},\ T\to\infty. 
\ee 
\end{mydef52} 

\begin{remark}
As a consequence, the Jacob's ladder can be viewed as complementary function to the function 
\bdis 
(1-c)\pi(T) 
\edis 
in the sense that 
\be \label{4.2} 
\vp_1(T)+(1-c)\pi(T)\sim T,\ T\to\infty. 
\ee 
\end{remark} 

\subsection{} 

Let (see \cite{3}, (1.11)) 
\be \label{4.3} 
\begin{split}
& y=\vp_1(t):\ \vp_1^0(t)=t,\ \vp_1^1(t)=\vp_1(t),\ \vp_1^2(t)=\vp_1(\vp_1(t)),\dots , \\ 
& \vp_1^k(t)=\vp_1(\vp_1^{k-1}(t)),\dots,\ t\in [T,T+U],\ T\geq T_0[\vp_1], 
\end{split}
\ee 
of course (see (\ref{4.1})) 
\be \label{4.4} 
T_0>\vp_1(T_0), 
\ee  
and the symbol $\vp_1^k(t)$ represents the $k$-th iteration of the Jacob's ladder. 

\begin{remark}
Let us remind that the functions 
\bdis 
\vp_1^k(t),\ k=2,3,\dots 
\edis  
are increasing since $\vp_1(t)$ is increasing. 
\end{remark} 

\subsection{} 

In the case 
\bdis 
t\mapsto \vp_1^k(t),\ t\in [T,T+U] 
\edis  
it follows from Lemma 2 that: 
\be \label{4.5} 
\vp_1^k(t)-\vp_1^{k+1}(t)\sim (1-c)\frac{\vp_1^k(t)}{\ln\vp_1^k(t)},\ k=0,1,\dots,n,\ t\to\infty, 
\ee  
where $n\in\mbb{N}$ is arbitrary and fixed one. Now formulae (\ref{4.5}) imply the following properties of the set 
\bdis 
\{\vp_1^k(t)\}_{k=0}^{n+1}\quad : 
\edis 

\begin{mydef53}
For 
\be\label{4.6} 
t\in [T,T+U],\ U=o\left(\frac{T}{\ln T}\right),\ T\to\infty 
\ee 
the following statements hold true: 
\be \label{4.7} 
t\sim \vp_1^1(t)\sim\vp_1^2(t)\sim \dots\sim\vp_1^{n+1}(t), 
\ee  
\be \label{4.8}
t>\vp_1^1(t)>\vp_1^2(t)>\dots>\vp_1^{n+1}(t), 
\ee 
\be \label{4.9} 
\vp_1^k(T)>(1-\epsilon)T,\ k=0,1,\dots,n+1,\ \epsilon>0, \ \epsilon\ \mbox{small and fixed}, 
\ee 
\be \label{4.10} 
\vp_1^k(T+U)-\vp_1^k(T)<\frac{1}{2n+5}\frac{T}{\ln T},\ k=1,\dots,n+1, 
\ee  
\be \label{4.11} 
\vp_1^k(T)-\vp_1^{k+1}(T+U)> 0.18 \times \frac{T}{\ln T},\ k=0,1,\dots,n. 
\ee 
\end{mydef53}

\subsection{} 

Further, we have introduced (see \cite{3}, (2.2)) the following set 
\be \label{4.12} 
D(T,U,n)=\bigcup_{k=0}^{n+1} [\vp_1^k(T),\vp_1^{k}(T+U)]. 
\ee  

\begin{remark}
We list here the properties of the set (\ref{4.12}): 
\begin{itemize}
	\item[(a)] It is disconnected set (see (\ref{4.11})) for every admissible $U$, (see (\ref{4.6})); 
	\item[(b)] Components of the set $D$ are distributed as follows: (see (\ref{4.11})) 
	\be \label{4.13} 
	\begin{split}
	& [\vp_1^{n+1}(T),\vp_1^{n+1}(T+U)]\prec [\vp_1^{n}(T),\vp_1^{n}(T+U)]\prec \dots \prec \\ 
	& \prec [\vp_1^{1}(T),\vp_1^{1}(T+U)]\prec [\vp_1^{0}(T),\vp_1^{0}(T+U)]=[T,T+U]. 
	\end{split}
	\ee 
\end{itemize}
\end{remark}

\begin{remark}
Asymptotic behaviour of the set $D$ is as follows: at $T\to\infty$ its components receding unboundedly each from other (see (\ref{4.11})) and all together recede to infinity. Hence at large $T$ the set (\ref{4.12}) behaves like one-dimensional Friedmann-Hubble expanding universe. 
\end{remark} 

\section{On the function $\tilde{Z}^2(t)$} 

Let us recall the following formula we have proved in \cite{1}: 
\be \label{5.1} 
Z^2(t)=\Phi'_\vp[\vp(t)]\frac{{\rm d}\vp(t)}{{\rm d}t},\ t\in [T,T+U],\ U=o\left(\frac{T}{\ln T}\right), 
\ee  
where 
\be \label{5.2} 
\begin{split}
& \Phi'_\vp[\vp]=\frac{2}{\vp^2}\int_0^{\mu[\vp]}t e^{-\frac{2}{\vp}t}Z^2(t){\rm d}t+ \\ 
& + Z^2\{\mu[\vp]\}e^{-\frac{2}{\vp}\mu[\vp]}\frac{{\rm d}\mu[\vp]}{{\rm d}\vp}, 
\end{split}
\ee 
(see \cite{1}, (3.5), (3.9)). Now we put (see (\ref{2.4}) and \cite{2}, (9.1)) 
\be \label{5.3}
\tilde{Z}^2(t)=\frac{{\rm d}\vp_1(t)}{{\rm d}t},\ t\geq T_0[\vp_1]. 
\ee 
In the next step we present just the result (see \cite{2}, Lemma 1, (7.7) -- (7.9), (9.2)): 

\begin{mydef54} 
If 
\be \label{5.4} 
\begin{split}
& \mu_a[\vp]=a\vp\ln\vp,\ a\in [7,8], \\ 
& t\in [T,T+U],\ U=o\left(\frac{T}{\ln T}\right), 
\end{split}
\ee  
then 
\be \label{5.5} 
\Phi'_\vp[\vp(t)]=\frac 12\left\{1+\mcal{O}\left(\frac{\ln\ln t}{\ln t}\right)\right\}\ln t, 
\ee  
i.e. (see (\ref{5.1}), (\ref{5.3})) 
\be \label{5.6} 
\tilde{Z}^2(t)=\frac{{\rm d}\vp_1(t)}{{\rm d}t}=\frac{1}{\omega(t)}\left|\zf\right|^2, 
\ee  
where 
\be \label{5.7} 
\omega(t)=2\Phi'_\vp[\vp(t)]=\left\{1+\mcal{O}\left(\frac{\ln\ln t}{\ln t}\right)\right\}\ln t,\ t\to\infty. 
\ee  
\end{mydef54} 

\begin{remark}
The segmant $[7,8]$ is sufficient one for our purposes since the continuum set $\mu_a[\vp]$ corresponds to this one (comp. (\ref{2.3}), (\ref{5.4})).  
\end{remark} 

\section{Reverse iterations}

\subsection{} 

Next, in our paper \cite{4} we have introduced the reverse iterations by means of the Jacob's ladder. First, we define the sequence 
\be \label{6.1} 
\{\overset{k}{T}\}_{k=0}^{k_0} 
\ee 
by the formula 
\be \label{6.2} 
\vp_1(\overset{k}{T})=\overset{k-1}{T},\ k=1,2,\dots,k_0,\ \overset{0}{T}=T,\ T\geq T_0[\vp_1] 
\ee 
(where $k_0\in\mbb{N}$ is an arbitrary and fixed numer) since the function $\vp_1(T),\ T\to\infty$ increases to $+\infty$. Further we have 
\be \label{6.3} 
\vp_1(\overset{k}{T})=\overset{k-1}{T} \ \Rightarrow \ \dots \ \Rightarrow \vp_1^k(\overset{k}{T})=T,\ k=1,\dots,k_0. 
\ee 
Since 
\be \label{6.4} 
\vp_1(\overset{1}{T})=T \ \Rightarrow \ \overset{1}{T}=\vp_1^{-1}(T), 
\ee  
and then we may use the inverse function $\vp_1^{-1}(T)$ to generate a reverse iterations. Namely we have: 
\be \label{6.5} 
\begin{split}
& \vp_1(\overset{2}{T})=\overset{1}{T} \ \Rightarrow \ \overset{2}{T}=\vp_1^{-1}(\overset{1}{T})=
\vp_1^{-1}(\vp_1^{-1}(T))=\vp_1^{-2}(T),\dots , \\ 
& \overset{k}{T}=\vp_1^{-k}(T),\ k=1,\dots,k_0, 
\end{split} 
\ee 
where the last row gives the $k$-th reverse iteration of the point $T=\overset{0}{T}$. Of course, we have 
\bdis 
\vp_1^k(\overset{k}{T})=\vp_1^k(\vp_1^{-k}(T))=T. 
\edis  

\subsection{} 

Now, the basic formula (\ref{4.1}) gives the following properties of the reverse iterations (see \cite{4}, (5.1) -- (5.13)). 

\begin{mydef55} 
If 
\be \label{6.6} 
U=o\left(\frac{T}{\ln T}\right),\ T\to\infty 
\ee  
then 
\be \label{6.7} 
\vp_1\{[\overset{k}{T},\overset{k}{\wideparen{T+U}}]\}=[\overset{k-1}{T},\overset{k-1}{\wideparen{T+U}}],\ [\overset{0}{T},\overset{0}{\wideparen{T+U}}]=[T,T+U], 
\ee 
\be \label{6.8} 
|[\overset{k}{T},\overset{k}{\wideparen{T+U}}]|=\overset{k}{\wideparen{T+U}}-\overset{k}{T}=o\left(\frac{T}{\ln T}\right), 
\ee  
\be \label{6.9} 
|[\overset{k-1}{\wideparen{T+U}},\overset{k}{T}]|=\overset{k}{T}-\overset{k-1}{\wideparen{T+U}}\sim (1-c)\frac{T}{\ln T}, 
\ee  
\be \label{6.10} 
[T,T+U]\prec [\overset{1}{T},\overset{1}{\wideparen{T+U}}]\prec \dots \prec [\overset{k}{T},\overset{k}{\wideparen{T+U}}],\ k=1,\dots,k_0, 
\ee 
(comp. Lemma 3 and (\ref{4.13})). 
\end{mydef55} 

From (\ref{6.6}) -- (\ref{6.10}) we obtain the following property of the Jacob's ladders. 

\begin{mydef82}
For every segment 
\bdis 
[T,T+U],\ U=o\left(\frac{T}{\ln T}\right),\ T\to\infty 
\edis  
there is the following class of disconnected sets (comp. (\ref{4.12})) 
\be \label{6.11} 
\Delta(T,U,k)=\bigcup_{r=0}^k [\overset{r}{T},\overset{r}{\wideparen{T+U}}],\ 1\leq k\leq k_0, 
\ee  
generated by the Jacob's ladder $\vp_1(T)$. 
\end{mydef82} 

\begin{remark}
Asymptotic behaviour of the set $\Delta$ is the same as behavior of the set (\ref{4.12}), i.e.  at $T\to\infty$ its components receding unboundedly each from other and all together recede to infinity. Hence at large $T$ the set (\ref{6.11}) behaves like one-dimensional Friedmann-Hubble expanding universe. 
\end{remark} 

\subsection{} 

Further, we have the following statement, see \cite{4}, (6.4): 

\begin{mydef56}
If 
\be \label{6.12} 
t\in [\vp_1^{-k}(T),\vp_1^{-k}(T+U)],\ k=1,\dots,k_0, 
\ee  
then (see (\ref{6.5})) 
\be \label{6.13} 
\vp_1^r(t)\in[\vp_1^{r-k}(T),\vp_1^{r-k}(T+U)],\ r=0,1,\dots,k, 
\ee  
i.e. 
\be \label{6.14} 
\begin{split}
& \vp_1^0(t)=t\in [\vp_1^{-k}(T),\vp_1^{-k}(T+U)]=[\overset{k}{T},\overset{k}{\wideparen{T+U}}], \\ 
& \vp_1^1(t)\in [\vp_1^{-k+1}(T),\vp_1^{-k+1}(T+U)]=[\overset{k-1}{T},\overset{k-1}{\wideparen{T+U}}], \\ 
& \vdots \\ 
& \vp_1^{k-1}(t)\in [\vp_1^{-1}(T),\vp_1^{-1}(T+U)]=[\overset{1}{T},\overset{1}{\wideparen{T+U}}], \\ 
& \vp_1^k(t)\in [\vp_1^{0}(T),\vp_1^{0}(T+U)]=[\overset{0}{T},\overset{0}{\wideparen{T+U}}]=[T,T+U]. 
\end{split}
\ee 
\end{mydef56}

\section{Main lemma and proof of Theorem 1}

\subsection{} 

In connection with direct and reverse iterations we have proved (see \cite{4}, (7.1), (7.2)) the following 

\begin{mydef57}
	If 
	\be \label{7.1}  
	U=o\left(\frac{T}{\ln T}\right),\ T\to\infty , 
	\ee 
	then for every Lebesgue-integrable function 
	\bdis    
	g(t),\ t\in [T,T+U] 
	\edis 
	the following holds true: 
	\be \label{7.2} 
	\int_T^{T+U} g(t){\rm d}t=\int_{\overset{k}{T}}^{\overset{k}{\wideparen{T+U}}}g[\vp_1^k(t)]\prod_{r=0}^{k-1}\tilde{Z}^2[\vp_1^r(t)]{\rm d}t,\ k=1,\dots,k_0. 
	\ee  
\end{mydef57}  

\begin{remark}
	We have obtained the case $k=1$: 
	\be \label{7.3} 
	\int_T^{T+U} g(t){\rm d}t=\int_{\overset{1}{T}}^{\overset{1}{\wideparen{T+U}}}g[\vp_1(t)]\tilde{Z}^2(t){\rm d}t 
	\ee 
	in our paper \cite{2}, (9.5). 
\end{remark} 

\subsection{} 

Now we proceed to the proof of our Theorem 1. 

\subsection*{Proof of Theorem 1} 

Since the system (\ref{1.3}) is fixed one then the corresponding $l$ is also fixed and consequently the condition (\ref{7.1}) 
\bdis 
l=o\left(\frac{T}{\ln T}\right),\ T\to\infty
\edis 
is fulfilled for all sufficiently big positive $T$. Now we have 
\bdis 
m\not=n:\ 0=\int_a^{a+2l}f_m(t)f_n(t){\rm d}t=\int_T^{T+2l}f_m(\tau-T+a)f_n(\tau-T+a){\rm d}\tau=
\edis 
and next, by our Lemma 7 for sufficiently big $T>0$, one obtains 
\bdis 
=\int_{\overset{p}{T}}^{\overset{p}{\wideparen{T+2l}}}f_m[\vp_1^p(\rho)-T+a]f_n[\vp_1^p(\rho)-T+a]\prod_{r=0}^{p-1}\tilde{Z}^2[\vp_1^r(\rho)]{\rm d}\rho= 
\edis 
and next, by simple substitution 
\bdis 
\rho=\frac{\overset{p}{\wideparen{T+2l}}-\overset{p}{T}}{2l}(t-a)+\overset{p}{T}, t\in [a,a+2l],\ \rho\in [\overset{p}{T},\overset{p}{\wideparen{T+2l}}]
\edis  
we obtain  
\be \label{7.4} 
\begin{split}
	& \frac{\overset{p}{\wideparen{T+2l}}-\overset{p}{T}}{2l}\int_{a}^{a+2l}f_m[\vp_1^p(\frac{\overset{p}{\wideparen{T+2l}}-\overset{p}{T}}{2l}(t-a)+\overset{p}{T})-T+a]\times \\ 
	& \times f_n[\vp_1^p(\frac{\overset{p}{\wideparen{T+2l}}-\overset{p}{T}}{2l}(t-a)+\overset{p}{T})-T+a]\times \\ 
	& \times \prod_{r=0}^{k-1}\tilde{Z}^2[\vp_1^p(\frac{\overset{p}{\wideparen{T+2l}}-\overset{p}{T}}{2l}(t-a)+\overset{p}{T})]{\rm d}t=
\end{split}
\ee 
and, in the next step, we finish with (see (\ref{1.5})) 
\be \label{7.6} 
=\frac{\overset{p}{\wideparen{T+2l}}-\overset{p}{T}}{2l}\int_a^{a+2l}f_m^p(t)f_n^p(t){\rm d}t \ \Rightarrow \ \boxed{\int_a^{a+2l}f_m^p(t)f_n^p(t){\rm d}t=0}. 
\ee  

\subsection{} 

We give at once the following. 

\begin{remark}
	If we use the last formula in (\ref{7.6}) as the origin of a new process (an analogue of this in the subsection 7.2), then we obtain $k^2$ new iterated $L_2$-orthogonal systems 
	\bdis 
	\{f_n^{p_1,p_2}(t)\}_{n=0}^\infty,\ t\in [a,a+2l],\ p_1,p_2=1,\dots,k 
	\edis  
	where 
	\bdis 
	\begin{split}
		& f_n^{p_1,p_2}(t)= \\ 
		& = f_n[\vp_1^{p_1}(\frac{\overset{p_1}{\wideparen{T+2l}}-\overset{p_1}{T}}{2l}(\vp_1^{p_2}(\frac{\overset{p_2}{\wideparen{T+2l}}-\overset{p_2}{T}}{2l}(t-a)+\overset{p_2}{T})-T)+\overset{p_1}{T})-T+a] \times \\ 
		& \times \prod_{r=0}^{p_1-1} \left|\tilde{Z}[\vp_1^{r}(\frac{\overset{p_1}{\wideparen{T+2l}}-\overset{p_1}{T}}{2l}(\vp_1^{p_2}(\frac{\overset{p_2}{\wideparen{T+2l}}-\overset{p_2}{T}}{2l}(t-a)+\overset{p_2}{T})-T)+\overset{p_1}{T})]\right| \times \\  
		& \times \prod_{r=0}^{p_2-1}
		\left|\tilde{Z}[\vp_1^r(\frac{\overset{p_2}{\wideparen{T+2l}}-\overset{p_2}{T}}{2l}(t-a)+\overset{p_2}{T})]\right|, 
	\end{split}
	\edis  
	and so on up to $k^l$ new iterated $L_2$-orthogonal systems 
	\bdis 
	\{f_n^{p_1,p_2,\dots,p_l}(t)\},\ t\in [a,a+2l],\ p_1,\dots,p_l=1,\dots,k
	\edis 
	for every fixed $l\in\mbb{N}$. 
\end{remark} 

\subsection{} 

Let us notice that the transformation (\ref{7.5}) 
\be \label{7.7} 
w=w(t)=\vp_1^p(\frac{\overset{p}{\wideparen{T+2l}}-\overset{p}{T}}{2l}(t-a))-T+a,\ t\in [a,a+2l]
\ee 
has the following properties: 
\begin{itemize}
	\item[(a)] by the subsection 6.1: 
	\bdis 
	\begin{split}
		& w(a)=\vp_1^p(\overset{p}{T})-T+a=T-T+a=a;\ \overset{p}{T}=\vp_1^{-p}(T), \\ 
		& w(a+2l)=\vp_1^p(\overset{p}{\wideparen{T+2l}})=a+2l, 
	\end{split}
	\edis 
	\item[(b)] since the function $\vp_1^p(u)$ is increasing one and 
	\bdis 
	u=\frac{\overset{p}{\wideparen{T+2l}}-\overset{p}{T}}{2l}(t-a)+\overset{p}{T},\ u\in [a,a+2l] 
	\edis 
	is evident then the composed function 
	\bdis 
	w(t),\ t\in [a,a+2l] 
	\edis  
	is increasing that is 
	\be \label{7.8} 
	w(t)\in [a,a+2l] . 
	\ee 
\end{itemize} 

\begin{remark}
	Consequently, it follows from (a) and (b) that by the one-to-one correspondence (\ref{7.7}) we have defined new automorphism on $[a,a+2l]$. i.e. the $k$, $k^2$, \dots , $k^l$ of new automorphisms for every fixed sufficiently big positive $T$. 
\end{remark}

I would like to thank Michal Demetrian for his moral support of my study of Jacob's ladders.

\end{document}